\documentclass[11pt]{amsart}
\usepackage{tabularx,booktabs,tikz}
\usepackage{caption}
\usepackage{amsmath}
\usepackage{amsfonts}

\usepackage{amscd}
\usepackage{amsthm}
\usepackage{amssymb} \usepackage{latexsym}
\usepackage{eufrak}
\usepackage{euscript}
\usepackage{epsfig}
\usepackage{graphics}
\usepackage{array}
\usepackage{enumerate}
\usepackage{dsfont}
\usepackage{color}
\usepackage{wasysym}
\usepackage{hyperref}
\usepackage{pdfsync}
\usepackage{tikz}

\newcommand{\bel}[1]{\begin{equation}\label{#1}}

\newcommand{\be}{\begin{equation}}

\newcommand{\ba}{\begin{eqnarray}}
\newcommand{\ea}{\end{eqnarray}}

\newcommand{\qe}{\end{equation}}
\newcommand{\rv}{{\mathbb R}^V}

\newcommand{\N}{{\mathbb N}}

\newcommand{\R}{{\mathbb R}}

\newcommand{\om}{\Omega}

\newcommand{\g}{\mathcal{G}}
\newcommand{\h}{\mathcal{H}}
\newcommand{\f}{\mathcal{F}}

\newcommand{\dg}{{\delta\g}}

\newcommand{\Hmm}[1]{\leavevmode{\marginpar{\tiny%
$\hbox to 0mm{\hspace*{-0.5mm}$\leftarrow$\hss}%
\vcenter{\vrule depth 0.1mm height 0.1mm width \the\marginparwidth}%
\hbox to
0mm{\hss$\rightarrow$\hspace*{-0.5mm}}$\\\relax\raggedright #1}}}

\newtheorem{theorem}{Theorem}[section]

\newtheorem{lemma}[theorem]{Lemma}
\newtheorem{corollary}[theorem]{Corollary}
\newtheorem{definition}[theorem]{Definition}

\newtheorem{remark}[theorem]{Remark}

\newtheorem{prop}[theorem]{Proposition}
\newtheorem{problem}[theorem]{Problem}

\newtheorem{example}[theorem]{Example}
\newtheorem{claim}[theorem]{Claim}

\newcommand{\tm}{\begin{theorem}}
\newcommand{\tmd}{\end{theorem}}
\newcommand{\co}{\begin{corollary}}
\newcommand{\cod}{\end{corollary}}
\newcommand{\prp}{\begin{prop}}
\newcommand{\prpd}{\end{prop}}
\newcommand{\pf}{\begin{proof}}
\newcommand{\pfd}{\end{proof}}
\newcommand{\rmk}{\begin{remark}}
\newcommand{\rmkd}{\end{remark}}
\newcommand{\ex}{\begin{example}}
\newcommand{\exd}{\end{example}}
\newcommand{\pr}{\begin{problem}}
\newcommand{\prd}{\end{problem}}
\newcommand{\df}{\begin{definition}}
\newcommand{\dfd}{\end{definition}}
\newcommand{\lm}{\begin{lemma}}
\newcommand{\lmd}{\end{lemma}}

\begin{document}

\title[Steklov flows on trees]{Steklov flows on trees and applications}

\author{Zunwu He}
\address{Zunwu He: School of Mathematics, South China University of Technology, 510641, Guangzhou, China}
\email{hzwmath789@scut.edu.cn}
%\address{Yohji Akama: Mathematical Institute, Graduate School of Science, Tohoku University,
%Sendai, 980-0845, Japan}
%\email{akama@math.tohoku.ac.jp}

\author{Bobo Hua}
\address{
Bobo Hua: School of Mathematical Sciences, LMNS, Fudan University, Shanghai 200433, China;  Shanghai Center for Mathematical Sciences, Jiangwan Campus, Fudan University, No. 2005 Songhu Road, Shanghai 200438, China.}
\email{bobohua@fudan.edu.cn}

%\author{Yanhui Su}
%\email{suyh@fzu.edu.cn}
%\address{Yanhui Su: College of Mathematics and Computer Science, Fuzhou University, Fuzhou 350116, China}

\begin{abstract}

%In this paper, we study the bounds for discrete Steklov eigenvalues on trees via geometric quantities. For a finite tree, we prove a sharp upper bound for the first nonzero Steklov eigenvalue by the reciprocal of the size of the boundary. This can be improved to the reciprocal of the size of the set of vertices for those trees with degree at least three for interior vertices.

%On the other hand, we obtain a sharp universal upper bound for the finite tree, in terms of the diameter of the tree.
%Moreover, we prove similar estimates for higher order Steklov eigenvalues.

We introduce the Steklov flows on finite trees, i.e. the flows (or currents) associated with the Steklov problem. By constructing appropriate Steklov flows, we prove the monotonicity and rigidity of the first nonzero Steklov eigenvalues on trees: for finite trees $\g_1$ and $\g_2,$ the first nonzero Steklov eigenvalue of $\g_1$ is greater than or equal to that of $\g_2$, provided that $\g_1$ is a subgraph of $\g_2.$ Moreover, we give the sufficient and necessary condition in which the equality holds.
\end{abstract}
\maketitle

%\tableofcontents
%\tableofcontents
%Mathematics Subject Classification 2010: 05C10, 31C05.

\par
\maketitle

\bigskip

%\begin{enumerate}
%\item all first names are omitted.
%\end{enumerate}

\section{introduction}

Given a $m$-dimensional ($m\geq 2$)
 compact orientable Riemannian manifold $(M,g)$ with smooth boundary $\partial M$, the Steklov problem on $(M,g)$ reads as
\begin{align*}
\left\{
\begin{array}{lr}
\Delta f(x)=0,\ \ x\in M,&
\\\frac{\partial f}{\partial n}(x)=\lambda f(x),\ x\in \partial M,&
\end{array}
\right.
\end{align*}
where $\Delta$ is the Laplace-Beltrami operator on $(M,g)$ and $\frac{\partial}{\partial n}$ is the outward normal derivative along $\partial M$. The spectrum of {the} Steklov problem on $(M,g)$ coincides with that of the following Dirichlet-to-Neumann operator \cite{Kuznetsov14},
\begin{align*}
\Gamma:&H^{\frac{1}{2}}(\partial M)\longrightarrow H^{-\frac{1}{2}}(\partial M)
\\&f\longmapsto\Gamma f:=\frac{\partial \hat f}{\partial n},
\end{align*}
where $\hat f$ is the harmonic extension to $M$ of $f.$
%\cite{Ne\v{c}as12}
It is well-known that the Dirichlet-to-Neumann operator is a first-order elliptic pseudo-differential operator \cite{Taylor11}, which is self-adjoint and non-negative. The spectrum of $\Gamma$ is discrete, and can be ordered as $$0=\lambda_1<\lambda_2\leq \lambda_3\leq\cdots\nearrow\infty,$$ where $\lambda_2$ is called the first (nonzero) Steklov eigenvalue. See \cite{Weinstock54,Brock01,Fraser11, Girouard12,Escobar97,Escobar99,Escobar00,Colbois11} for more results of Steklov eigenvalues on Riemannian manifolds.

%{%For the Steklov eigenvalues on Riemannian manifolds,

In order to detect spectral properties of Riemannian manifolds, Colbois et al. investigated the Steklov problem on some discretizations of manifolds \cite{Colbois18}; see Section~\ref{sec:pre} for the definition of the Steklov problem on graphs. The second author, Huang and Wang \cite{Hua17}, and Hassannezhad and Miclo \cite{Miclo2017}, studied the first Steklov eigenvalue on graphs using isoperimetric constants independently. Some lower bound estimates of the first Steklov eigenvalue on graphs were proved in \cite{Perrin19,ShiYu20a}. For subgraphs in Cayley graphs of discrete groups of polynomial growth, the upper bound estimates were proved in \cite{Perrin20}; see also \cite{Han19}. In our previous paper \cite{HH20}, we obtained various upper bounds of Steklov eigenvalues on finite trees. Note that infinite trees are regarded as discrete counterparts of Hadamard manifolds. See \cite{HuaHuangWang18,ShiYu19a,ShiYu19b,ShiYu20b} for other developments on Steklov eigenvalues of graphs.

In this paper, we study the Steklov eigenvalues for the Steklov problem on finite trees.
For investigating the relation of the first Steklov eigenvalues of finite trees, we introduce a $\lambda$-flow to some interior vetex or some boundary vertex on a finite tree $\g$ (see Definition \ref{lambda-flow}), which can be regarded as a generalization of Steklov eigenfunction on $\g$ associated with eigenvalue $\lambda$. For a finite tree $\g=(V,E),$ we denote by $\dg$ the boundary of $\g$, i.e. the set of pending vertices, and by $\om:=V-\dg$ the set of interior vertices of $G.$
\df\label{lambda-flow}
For a finite tree $\g=(V,E)$ with boundary $\dg$, if there is a nonzero function $f_\lambda\in \rv$, some $\lambda\geq 0$ and $x\in V$ such that

%\item $\nabla_{(y,z)}f=f(y)-f(z)>0$ for every $(y,z)\in E$ with $d(x,y)>d(x,z)$
\begin{equation}
\left\{
\begin{array}{lr}
 \Delta f_\lambda(w)=0, \quad \quad w\in \om-\{x\}, &\\
 \frac{\partial f_\lambda}{\partial n}(w)=\lambda f(w), \ w\in \dg-\{x\},&
\end{array}
\right.
\end{equation}
then we say that $f_\lambda$ is a $\lambda$-flow to $x$ on $\g$, or simply a $\lambda$-flow to $x.$ In the above, $\frac{\partial}{\partial n}$ is the discrete version of outward normal derivative; see Section~\ref{sec:pre}.
\dfd
\rmk
\begin{enumerate}
\item For a $\lambda$-flow $f_\lambda$ to $x$ on $\g,$ if $\sum\limits_{z\in\dg}f_\lambda(z)=0$ in case of $x\in\dg$, or $\Delta  f_\lambda(x)=0$ in case of $x\in \om$, then $f_\lambda$ is a Steklov eigenfunction on $\g$ associated with eigenvalue $\lambda$.
\item There are fewer constraint equations for a $\lambda$-flow than those for a Steklov eigenfunction. In fact, the $\lambda$-flow $f_\lambda$ has one degree of freedom. This guarantees the existence of $\lambda$-flows on $\g$ for $\lambda$ in some interval, and these $\lambda$-flows are continuous in $\lambda;$ see Lemma \ref{exixt posi}.
\item The flows (or currents) are well-studied in the literature of electrical networks; see e.g. \cite[Definition~2.1]{Barlowbook}. For $f_\lambda$ defined above, we regard it as the potential in an electrical network, i.e. the voltage. In this paper, we would rather consider the flow induced by the potential in the proofs for convenience, and hence call $f_\lambda$ the $\lambda$-flow.
\end{enumerate}% as $\lambda$ is in some interval.
\rmkd

%\rmk
%Note that $cf_\lambda$ is also a $\lambda$-flow to $x$ on $\g$, for any $c\neq 0$. If $x$ is a interior vertex of $V$, i.e $x\in \om$, then $\lambda$ is a Steklov eigenvalue of the finite tree $\g$ and the $\lambda$-flow $f_\lambda$ is a Steklov eigenfunction for $\lambda$. But if $x\in \dg$, $\lambda$ is not a Steklov eigenvalue of $\g$ and $f_\lambda$ is not a eigenfunction in general.
%\rmkd

The main result is formulated as follows. For two graphs $\g_1$ and $\g_2,$ we say that $\g_1$ is a subgraph of $\g_2$ if there is an injective graph homomorphism from $\g_1$ to $\g_2.$
\tm\label{main1}
%\red{If two finite trees $\g_1=(V_1,E_1),\g_2=(V_2,E_2)$ satisfy $V_1\subset V_2,E_1\subset E_2$, then}

Given two finite trees $\g_1=(V_1,E_1)$ and $\g_2=(V_2,E_2),$ if $\g_1$ is a subgraph of $\g_2,$ then

\begin{align}
\lambda_2(\g_2)\leq\lambda_2(\g_1),
%\blue{where $\lambda_2$ is the first Steklov eigenvalue.}
\end{align}where $\lambda_2$ is the first Steklov eigenvalue.

Furthermore, the equality holds if and only if $\g_1=\g_2$ or there exist $x$, its neighbors $x_j$ in $V_1$ and its neighbors $y_s$ in $V_2$ such that $\sigma(\g_{1,1},x)=\sigma(\g_{1,2},x)=\min\limits_{1\leq j\leq m_1}\sigma(\g_{1,j},x)=\sigma(\g_{2,1},x)=\sigma(\g_{2,2},x)=\min\limits_{1\leq s\leq m_2}\sigma(\g_{2,s},x)$, where $m_1$ ($m_2$,resp.) is the degree of $x$ in $\g_1$ ($\g_2$, resp.), $\h_{1,j}$ ($\h_{2,s}$, resp.) is the branch from $(x_j,x)$ ($(y_s,x)$,resp.) in $\g_1$ ($\g_2$, resp.), and $\g_{1,j}:=\h_{1,j}(x_j,x)=(V_{1,j},E_{1,j})$ ($\g_{2,s}:=\h_{2,s}(y_s,x)=(V_{2,s},E_{2,s})$, resp.) with boundary $\dg_{1,j}$ ($\dg_{2,s}$, resp.).
%Let $\hat{\mathcal{\g}}:=(\g)^2_x=(\hat V,\hat E)$ (see Definition \ref{wedge sum}) with boundary $\delta\hat{\mathcal{\g}}$ and $g_2\in \rhv$ be any Steklov eigenfunction associated to eigenvalue $\lambda_2(\hat{\mathcal{\g}})$, then $g_2(x)=0$.
%\lambda_2(\mathcal{\g})
%Moreover, we have $\min\limits_{1\leq j\leq m}\sigma(\g_j,x)=\lambda_2(\hat{\mathcal{\g}})$}
\tmd

\begin{figure}[htbp]
 \begin{center}
   \begin{tikzpicture}
    \node at (0,0){\includegraphics[width=0.5\linewidth]{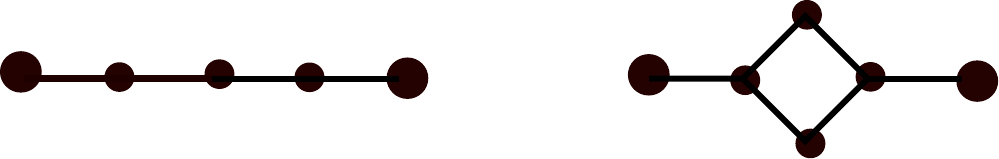}};
   % \node at (0, -.8){\Large $e$};
    \node at (-1.8,   -1){\Large $\g_1$};
    \node at (2,   -1){\Large $\g_2$};
  %  \node at (0,  1.3){\Large $\sigma_1$};
  %  \node at (0,  -1.7){\Large $\sigma_2$};
   % \node at (-3.4,  1.5){\Large $z$};
   \end{tikzpicture}
  \caption{The counterexample for monotonicity of the Steklov eigenvalues on general graphs.}\label{fig0}
 \end{center}
\end{figure}

\rmk
\begin{enumerate}
\item This result fails for general finite graphs; see Figure \ref{fig0}. Here $\g_1$ is a subgraph of $\g_2,$ but $\lambda_2(\g_1)=\dfrac{1}{2}<\lambda_2(\g_2)=\dfrac{2}{3}$.
\item The monotonicity for eigenvalues is important in the spectral theory. It is well-known that Dirichlet eigenvalues for subgraphs are monotone with respect to the inclusion of subgraphs, \cite[Theorem~2.3]{Friedman93}; see interlacing inequalities for general setting \cite{Haemers95}. For the class of finite trees, we prove the monotonicity of the first Steklov eigenvalue, which is possibly a consequence by the combinatorics of trees.
\item We appreciate that some referees could provide an interesting simplified proof of the part for the monotonicity. However, the rigidity part is quite subtle. We first give a complicate proof for the monotonicity, and then use the arguments to show the rigidity of the result. This is the novel contribution of the paper.
\end{enumerate}
%Direct computation gives .
\rmkd

One of the motivations of introducing $\lambda$-flows on the branches in a tree and associated $\sigma$-eigenvalues (see Definition \ref{sigma2}) is Theorem \ref{sigmasteklov}. It relates $\sigma$-eigenvalues and the first nonzero Steklov eigenvalue.

We sketch the proof strategy for Theorem \ref{main1}. Since one graph is a subgraph of the other, it suffices to consider the case of that $\g^\prime$ is obtained from $\g$ by adding one vertex $y$ and {directed edges ${(x,y),(y,x)}$.} The basic idea is that we construct a $\lambda$-flow $f_\lambda$ to $y$ along a path on $\g^\prime,$ which can be descended to a continuous ``homotopy mapping'' on $\g\times [0,\lambda_2(\g)]$ between $f_0|_{\g}$ (the eigenfunction on $\g$ associated to eigenvalue $0$) and $f_{\lambda_2(\g)}|_{\g}$ (the eigenfunction on $\g$ associated to eigenvalue $\lambda_2(\g)$). This can be used to deduce that $\sum\limits_{z\in\dg^\prime}f_{\hat\lambda}(z)=0$ with $0<\hat\lambda\leq\lambda_2(\g)$ and hence $f_{\hat\lambda}$ must be a Steklov eigenfunction on $\g^\prime$ associated with eigenvalue $\hat\lambda$.

In the following, we turn to the crucial construction of the $\lambda$-flow.
Lemma \ref{exixt posi} and Lemma \ref{subtree} are two key ingredients to construct the above $\lambda$-flow $f_\lambda$. Lemma \ref{subtree} implies that $f_\lambda$ must take the same sign on the relative boundaries of some ``sufficiently large'' branches, as $\lambda$ is small.
Let $z$ be a boundary vertex. Lemma \ref{exixt posi} shows the existence, uniqueness (up to scaling) and continuity of some $\lambda$-flow $g_\lambda$ to $z$ on $\g$ for sufficiently small $\lambda$, such that $g_\lambda$ has the same sign on the relative boundary vertices except $z$. These $\lambda$-flows satisfying Lemma \ref{exixt posi} are the fundamental blocks to construct the previous $\lambda$-flow to $y$ on $\g^\prime$.

The above $\lambda$-flow $f_\lambda$ to $y$ on $\g^\prime$ is constructed inductively by adding some $\lambda$-flows satisfying Lemma \ref{exixt posi} on subtrees of $\g^\prime$. In order to ensure that $f_{\lambda_2(\g)}|_{\g}$ is a Steklov eigenfunction on $\g$ associated with $\lambda_2(\g)$, Lemma \ref{subtree} indicates that we must find a ``deepest'' edge $(u,v)$ in $\g$ with two subtrees $\h_1,\h_2$ obtained by removing $(u,v)$, such that $f_\lambda|_{\g}$ takes positive values on $\h_1$ containing $u$ and non-positive values on $\h_2$ containing $v$.

\begin{figure}[htbp]
 \begin{center}
   \begin{tikzpicture}
    \node at (0,0){\includegraphics[width=0.5\linewidth]{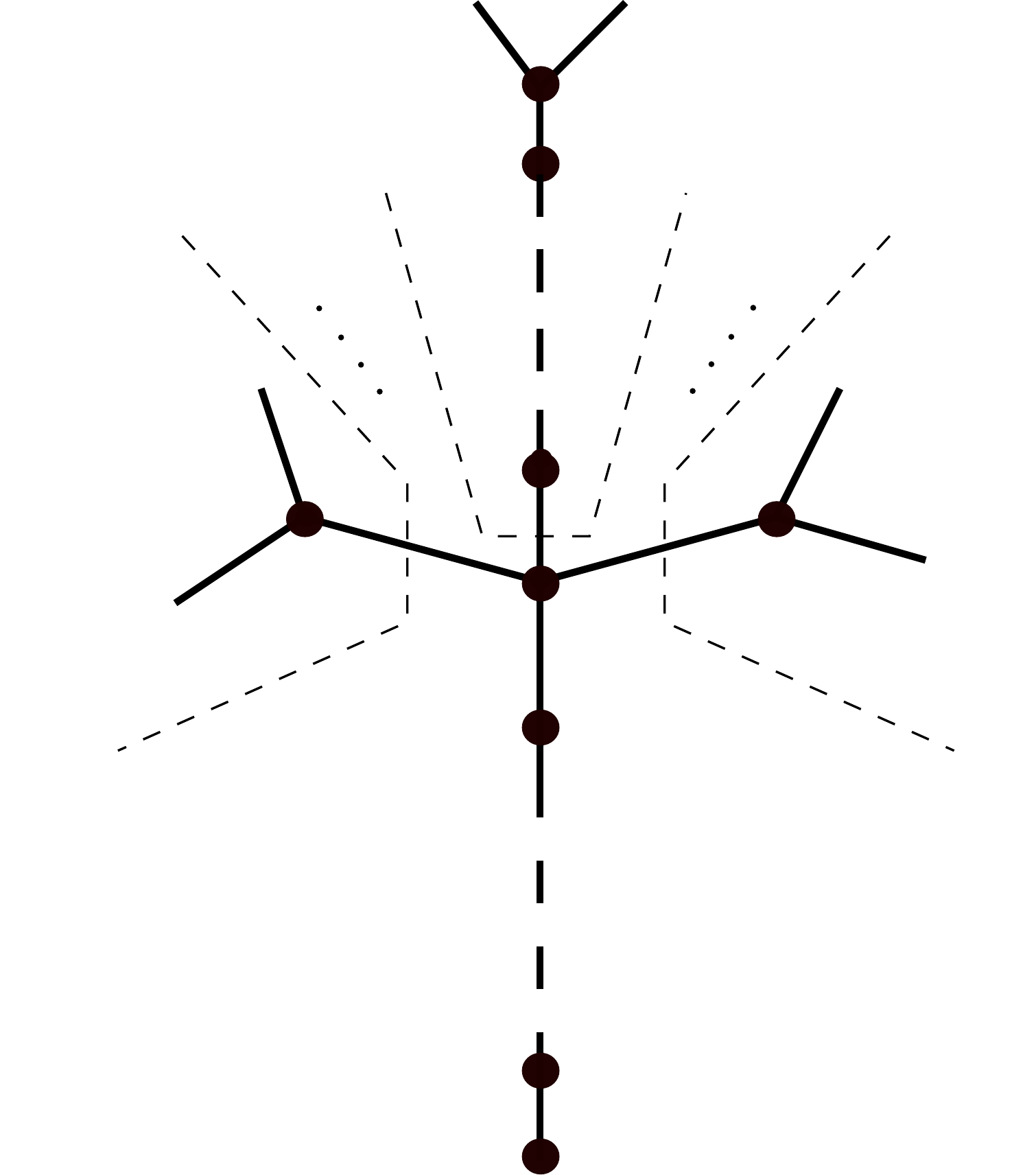}};
   % \node at (0, -.8){\Large $e$};
    \node at (-2.8,   -0.7){\Large $\h_{i,1}$};
    \node at (2.6,   -0.6){\Large $\h_{i,j}$};
    \node at (-0.3,   2){\Large $\g_i$};
     \node at (-0.2,   0.9){ $v_{i-1}$};
    \node at (-0.2,   -0.2){ $v_i$};
    \node at (-0.4,   -1){ $v_{i+1}$};
    \node at (-1.7,   0.5){$v_{i,1}$};
     \node at (2,   0.5){ $v_{i,j}$};
    %  \node at (0,  1.3){\Large $\sigma_1$};
  %  \node at (0,  -1.7){\Large $\sigma_2$};
   % \node at (-3.4,  1.5){\Large $z$};
   \end{tikzpicture}
  \caption{}\label{fig1}
 \end{center}
\end{figure}

This provides the first step of construction of the above $f_\lambda$. We may assume $\h_2$ contains $x$ and there is a path $\alpha:=v_0=v\sim v_1\sim v_2\sim\cdots\sim v_m=x\sim y$ connecting $v$ and $y$ in $\g^\prime$. By induction and Lemma \ref{exixt posi}, it admits a $\lambda$-flow $f_{\lambda,i,j}$ ($f_{\lambda,i}$ resp.) to $v_i$ on each subtree $\h_{i,j}(v_{i,j},v_i)$ ($\g_i(v_{i-1},v_i)$ resp.) with $d(v_{i,j},v_i)=1$, and $1\leq j\leq \deg(v_i)-2$ for some $1\leq i\leq m$, where $\h_{i,j}$ ($\g_i$ resp.) is the branch from $(v_{i,j},v_i)$ ($(v_{i-1},v_i)$ resp.) in $\g$ (some notations see Definition \ref{branch}); see Figure \ref{fig1}. Since $(u,v)$ is ``deepest'', the previous flow $f_{\lambda,i,j}$ takes nonzero value at $v_i$ as $0\leq\lambda\leq\lambda_2(\g)$. This implies $f_{\lambda,i,j}$ can be used to assemble a $\lambda$-flow on a larger subtree. We go forth inductively to construct the $\lambda$-flow $f_\lambda$ along the path $\alpha$, and finish the construction.

In our previous paper \cite{HH20}, we prove the upper bound estimate of the first {nonzero} Steklov eigenvalue using the diameter of a finite tree; see \cite[Theorem~1.4]{HH20}. As an application of the main result, we characterize the equality case for the upper bound estimate; see Theorem~\ref{optimal diam}.
Moreover, we prove a lower bound on the first Steklov eigenvalue of finite trees, in terms of the bounds for the vertex degree and the diameter; see Theorem~\ref{large tree}, and also characterize the equality case. These yield rigidity results for the estimates of the first Steklov eigenvalue via the diameter.

At the end of the introduction, we propose two open problems on the monotonicity of higher order Steklov eigenvalue of some finite trees and graphs.
\pr\label{prob1}
Let $\g_1,\g_2$ be two finite trees with boundaries $\dg_1,\dg_2$ such that $\g_1$ is a subgraph of $\g_2$. Does it hold that $\lambda_k(\g_1)\geq \lambda_k(\g_2)$ for $3\leq k\leq |\dg_1|$, where $\lambda_k$ is the $k$-th Steklov eigenvalue?
\prd

\pr\label{prob2}
Let $\g_1,\g_2$ be two finite graphs with boundaries $\dg_1,\dg_2$ such that $\g_1$ is a subgraph of $\g_2$. If $i_*:\pi_1(\g_1)\longrightarrow\pi_1(\g_2)$ is a group isomorphism, where the map $i$ is the canonical inclusion from $\g_1$ to $\g_2$ and $\pi_1$ denotes the fundamental group, does it hold that $\lambda_k(\g_1)\geq \lambda_k(\g_2)$ for $2\leq k\leq |\dg_1|$, where $\lambda_k$ is the $k$-th Steklov eigenvalue?
\prd

\rmk
 Problem \ref{prob2} is equivalent to the case that $\g_1=(V_1,E_1),\g_2=(V_2,E_2)$ with $V_2=V_1\cup\{y\},E_2=E_1\cup\{(x,y)\}$, where $x\in V_1,y\notin V_1$ and $(x,y)\in E_2$, $E_1,E_2$ are sets of undirected edges.

The assumption on $i_*$ can not be removed; see the counterexample in Figure \ref{fig0}.
\rmkd

% which is a rigidity result for such estimate.
%\red{We also propose a explicit lower bound on the first Steklov eigenvalue of finite trees, in terms of the bounds of the degree and the diameter. The realization of the low bound is described fully. The readers are referred to see Theorem \ref{large tree} for the precious statement.}

The paper is organized as follows: In next section, we recall some basic facts on graphs.
Section~\ref{sec:steklov} is devoted to the Steklov flow. In Section~\ref{sec:proof}, we prove the main result, Theorem \ref{main1}. The last section contains some applications of the main result.

\textbf{Added in the proof:}  After the submission of our paper, Yu and Yu provide positive answers to Problem \ref{prob1} and Problem \ref{prob2} in a recent preprint \cite{Yu21}.
%\co%Theorem \ref{comparison1}
%If two finite trees $\g_1=(V_1,E_1)$, $\g_2=(V_2,E_2)$ with $x_1\in V_1,x_2\in V_2$ such that $\lambda_2(\g_1)^2_{x_1}=\lambda_2(\g_2)^2_{x_2}$. Then we have

%\begin{align}
%\lambda_2(\g_1\cup\g_2/x_1\sim x_2)=\lambda_2((\g_1)^2_{x_1}).
%\end{align}
%\cod

%\co%Theorem \ref{main1} and Theorem \ref{comparison1}Given any finite tree $\h=(V_1,E_1)$ and $x\in V_1$, for every $k\geq 3$ we have
%\begin{align}
%\lambda_2((\h)^2_x)=\lambda_2((\h)^k_x).
%\end{align}
%\cod
\textbf{Acknowledgements.} B.H. is supported by NSFC, grants no.11831004 and no. 11926313.

%\textbf{Data Available Statement.} The data used to support the findings of this study are available from the corresponding author upon request.

\section{preliminaries}\label{sec:pre}
First, we introduce some basic concepts. Let $(V,E)$ be a graph with the set of vertices $V$ and the set of directed edges $E.$ In this paper, we denote by $|\cdot|$ the cardinality of a set.
\df\label{wedge sum}
Let $\g_i=(V_i,E_i),$ $i=1,2,$ be graphs and $x_i\in V_i.$  We say that the graph $\g_1\bigsqcup\g_2/x_1\backsim x_2:=(V,E)$ is the wedge sum of $\g_1,\g_2$ at $x_1,x_2$ if $V=V_1\bigsqcup V_2/{x_1\backsim x_2}$ and $E=E_1\bigsqcup E_2$. Here $\bigsqcup$ means disjoint union, $x_1\backsim x_2$ means the identification of two vertices.

In particular, we denote by $(\g)^2_{x}=(\hat V,\hat E)$ the wedge sum of $\g,\g'$ at $x,x'$ with $\g=(V,E),\g'=(V',E')$, $\hat V=V\bigsqcup V'/x\sim x',\hat E=E\bigsqcup E'$, where there is a map $\phi$ such that $x\in V\rightarrow \phi(x)=x'\in V'$ and $(x,y)\in E\rightarrow \phi((x,y))=(x',y')\in E'$, and $\phi$ is a graph isomorphism.
\dfd
%if there is a graph isomorphism $\phi:\g_1\longrightarrow\g_2$ such that $\phi(x_1)=x_2$,

\rmk\label{double}
We does not distinguish two graphs up to graph isomorphism. We regard $\g$ as a subgraph of $(\g)^2_{x}$ with $V\subset\hat V$, $E\subset\hat E$.
%Given a graph $\g=(V,E)$ with $x\in V$, since $id:\g\longrightarrow\g$ and $id(x)=x$, $(\g)^2_x$
\rmkd
From now on, we always consider finite trees $\mathcal{G}=(V,E)$, and $E$ is a set of directed edges. Denote by $(\mathcal{G},\dg)$ the graph $\mathcal{G}=(V,E)$ with the boundary $\dg$ and the interior $\Omega:=V\setminus\dg.$

For a subtree $\g_1=(V_1,E_1)$ in $\g$, we define the relative boundary $\delta(\g_1,\g)$ (boundary $\delta\g_1$, resp.) of $\g_1$ to be $\dg\cap V_1$ (all the vertices in $V_1$ of degree one in $\g_1$, resp.).

Note that $\Omega$ is connected and there is no edge connecting two boundary vertices if $|V|\geq 3$. For any subset $S\subset V,$ we denote by $\mathbb{R}^{S}$ the vector space of all real functions on $S$ over $\mathbb{R}.$ It can be equipped with $\ell^2$-inner product $(\cdot,\cdot)_S:$ for any $f,g\in \mathbb{R}^{S}$, set $(f,g)_S:=\sum\limits_{x\in S}f(x)g(x).$ Then $(\mathbb{R}^{S},(\cdot,\cdot)_S)$ is a Hilbert space. For $f\in \mathbb{R}^V$, one can define the Laplace operate $\Delta$ on $\mathbb{R}^{V}$ such that
\begin{align}\label{laplace}
(\Delta f)(x):=\sum\limits_{y\in V:y\sim x}(f(x)-f(y)),
\end{align} where $y\sim x$ means that they are neighbors.
For convenience, we set $\nabla_{(u,v)}f:=f(u)-f(v)$ for any $f\in\rv$.

The outward normal derivative operator is defined as
\begin{align*}
\frac{\partial}{\partial n}:&\mathbb{R}^V\longrightarrow \mathbb{R}^{\dg}
\\&f\longmapsto\frac{\partial f}{\partial n},
\end{align*}
where $\frac{\partial f}{\partial n}(x)=f(x)-f(y)$ for any $x\in \dg$ and $y\sim x$. Since there is no edge connecting two boundary vertices, $\frac{\partial f}{\partial n}(x)=(\Delta f)(x)$ for any $x\in \dg.$

We introduce the Steklov problem on the pair $(\mathcal{G},\dg)$. For any nonzero function $f\in \mathbb{R}^V$ and some $\lambda\in \mathbb{R}$, the following equations hold,
\begin{equation}\label{Steklov 2}
\left\{
\begin{array}{ll}
\Delta f(x)=0, &x\in \Omega,\\
\dfrac{\partial f}{\partial n}(x)=\lambda f(x),&x\in \dg.
\end{array}
\right.
\end{equation}
The above $\lambda$ is called the Steklov eigenvalue of the graph $\g$ with boundary $\dg$, and $f$ is a Steklov eigenfunction associated to eigenvalue $\lambda$.

In analogy to the Riemannian case, one can define the Dirichlet-to-Neumann operator in the discrete setting as
\begin{align}
\Lambda:&\mathbb{R}^{\dg}\longrightarrow \mathbb{R}^{\dg}\nonumber
\\&f\longmapsto \Lambda f:=\frac{\partial \hat f}{\partial n},\label{dtn}
\end{align}
where $\hat f$ is the harmonic extension of $f$. Namely, $\hat f$ satisfies
\begin{equation}\label{dtn2}
\left\{
\begin{array}{ll}
\Delta \hat f(x)=0, & x\in \Omega ,\\
\hat f(x)=f(x),  &x\in \dg.
\end{array}
\right.
\end{equation}
It is well known that the Steklov eigenvalues in (\ref{Steklov 2}) are exactly the Dirichlet-to-Neumann eigenvalues in (\ref{dtn}), and the Steklov eigenfunctions in (\ref{Steklov 2}) are  the harmonic extensions of the corresponding eigenfunctions in (\ref{dtn}). The Dirichlet-to-Neumann operator is non-negative and self-adjoint. Since $V$ is finite, there are $|\dg|$ Steklov eigenvalues. We may arrange the Steklov eigenvalues in the following way:
\begin{align}\label{steklov eigen}
0=\lambda_1\leq \lambda_2\leq \cdots\leq \lambda_{|\dg|}.
\end{align}

In fact, since $\Omega$ is connected, $\lambda_2>0$ and $\lambda_{|\dg|}\leq 1$ \cite{Hua17,HuaHuangWang18}. Moreover, $(f_i,f_j)_\dg=0$ if $\lambda_i\neq \lambda_j$, where $f_i$ and $f_j$ are eigenfunctions with eigenvalues $\lambda_i,\lambda_j$ respectively. Note that constant functions are the eigenfunctions of $\lambda_1=0$. The space $\mathbb{R}^{\dg}$ has an orthonormal basis consisting of the Dirichlet-to-Neumann eigenfunctions.

%For $f\in $
For any $0\neq f\in \rv$, recall that $E$ is directed edge set, we define the Rayleigh quotient
\begin{align}\label{ray}
R(f)=\dfrac{\sum\limits_{(x,y)\in E}(f(x)-f(y))^2}{2\sum\limits_{x\in \dg}f^2(x)}.
\end{align}
The variational principles for $\lambda_k$ read as
\begin{align}\label{lamk}
&\lambda_k=\min\limits_{W\subset \mathbb{R}^V,\ dimW=k}\max_{0\neq f\in W}R(f)
\\&\lambda_k=\min_{\substack{W\subset \mathbb{R}^V,dimW=k-1\\W\perp 1_{\dg}}}\max\limits_{0\neq f\in W}
R(f)\label{var k},
\end{align}
where $1_{\dg}$ is the characteristic function on $\dg$, i.e. $f(x)=1$ if $x\in \dg$ and $f(x)=0,$ otherwise, and $W\perp 1_{\dg}$ means that any function in $W$ is orthogonal to $1_{\dg}.$

\section{The Steklov flows}\label{sec:steklov}
In this section, we aim to introduce the Steklov flows, i.e. $\lambda$-flows, and study their properties.

\df\label{branch}
Given a finite tree $\g=(V,E)$ with $x,y\in V,(x,y)\in E$, we call the subtree $\h$ is the branch from $(x,y)$ if $\h$ is the connected component containing $x$ as one removes the edge $(x,y)$ in $\g$. Denote by $\h(x,y)$ the subtree generated by $\h\cup\{y\}$ in $\g$, i.e. the vertices of the subtree are exactly those of $\h$ and $\{y\}$.
\dfd
\rmk
Note that the other connected component which contains $y$ is the branch from $(y,x)$.
\rmkd

%\red{?????}\df\label{lambda-flow}
%For a finite tree $\g=(V,E)$ with boundary $\dg$, if there is a nonzero function $f_\lambda\in \rv$, some number $\lambda\geq 0$ and $x\in V$ such that

%\item $\nabla_{(y,z)}f=f(y)-f(z)>0$ for every $(y,z)\in E$ with $d(x,y)>d(x,z)$
%\begin{equation}
%\left\{
%\begin{array}{lr}
% \Delta f_\lambda(w)=0, \quad \quad w\in \om-\{x\} &\\
 %\frac{\partial f_\lambda}{\partial n}(w)=\lambda f(w), \ w\in \dg-\{x\}&
%\end{array}
%\right.
%\end{equation}

%We say $f_\lambda$ is a $\lambda$-flow to $x$ on $\g$, or simply a $\lambda$-flow to $x$ if there are no confusions from the context.
%\dfd
%\rmk
%Note that $cf_\lambda$ is also a $\lambda$-flow to $x$ on $\g$, for any $c\neq 0$.

%If $\sum\limits_{z\in\dg}f_\lambda(z)=0$ as $x\in\dg$, or $\Delta  f_\lambda(x)=0$ as $x\in V-\dg$, then $f_\lambda$ is a Steklov eigenfunction on $\g$ associated with eigenvalue $\lambda$.

%\rmkd

We give some useful properties for $\lambda$-flows.
\lm\label{scale}
For a finite tree $\g=(V,E)$ with boundary $\dg$ containing $x$, assume that $f_\lambda\in \rv$ is a $\lambda$-flow to $x$ for some $\lambda\geq 0$ with $f_\lambda(x)=0$.
Then for any $c\neq 0$, we have
$$\lambda=R(f_\lambda)=R(cf_\lambda).$$
\lmd
\pf
Note that $cf_\lambda$ is also a $\lambda$-flow to $x$. Since $\g$ is finite and $E$ is directed edge set, it always holds that for any $g\in \rv,$
\begin{align}\label{green formula}
\frac12\sum\limits_{(u,v)\in E}(g(u)-g(v))^2=(\Delta g,g)_V.
\end{align}
Hence by (\ref{ray}), (\ref{green formula}) and Definition \ref{lambda-flow}, we have
\begin{align*}
R(cf_\lambda)&=R(f_\lambda)=\dfrac{(\Delta f_\lambda,f_\lambda)_V}{(f_\lambda,f_\lambda)_\dg}
=\dfrac{(\Delta f_\lambda,f_\lambda)_\dg}{(f_\lambda,f_\lambda)_\dg}
\\&=\dfrac{(\Delta f_\lambda,f_\lambda)_{\dg-\{x\}}}{(f_\lambda,f_\lambda)_{\dg-\{x\}}}
=\lambda.
\end{align*}
\pfd

\lm\label{edge flow}
For a finite tree $\g=(V,E)$ with boundary $\dg$ containing $x$, assume that $f_\lambda\in \rv$ is a $\lambda$-flow to $x$ for some $\lambda\geq 0$ and $(u,v)\in E$ with $d(x,u)>d(x,v)$. Then
$$\nabla_{(u,v)} f_\lambda=\sum\limits_{z\in \delta\h(u,v)-\{v\}}f_\lambda(z)\lambda,$$
where $\h$ is the branch from $(u,v)$, $\delta\h(u,v)$ is the boundary of $\h(u,v)$.
\lmd
\pf
Assume that $v,u_j$ are all neighbours of $u$ for $1\leq j\leq m-1$ with $m=\deg(u)$. We argue by induction on $|V|$. The case of $|V|=2$ is trivial.

By induction on $\h(u_j,u)$, we have $\nabla_{(u_j,u)} f_\lambda=\sum\limits_{z\in \delta\h(u_j,u)-\{u\}}f_\lambda(z)\lambda.$ Since $\Delta f_\lambda(u)=0$ and $\delta\h(u,v)-\{v\}$ is the union of $\delta\h(u_j,u)-\{u\}$ for $1\leq j\leq m-1$, we deduce that $$\nabla_{(u,v)} f_\lambda=\sum\limits_{j=1}^{m-1}\nabla_{(u_j,u)}f_\lambda=\sum\limits_{z\in \delta\h(u,v)-\{v\}}f_\lambda(z)\lambda.$$

\pfd

\df\label{sigma2}
For a finite tree $\g=(V,E)$ with boundary $\dg$ containing $x$, if there exists a $\lambda$-flow to $x$, $f_\lambda\in \rv$ with $\lambda\geq 0$ such that
\begin{enumerate}[(i)]
\item $f_\lambda(x)=0,$ \label{sigma2,1}
\item $\nabla_{(y,z)} f_\lambda=f_\lambda(y)-f_\lambda(z)>0$ if $(y,z)\in E$ and $d(x,y)>d(x,z)$.\label{sigma2,2}
\end{enumerate}
Then we denote by $\Sigma(\g,x)$ $(\hat\Sigma(\g,x), resp.)$ the set of such $\lambda$'s satisfying (\ref{sigma2,1}), (\ref{sigma2,2}) ((\ref{sigma2,1}), resp.) and $$\sigma(\g,x):=\inf\limits_{\lambda\in \Sigma(\g,x)}\lambda\quad \left(\hat\sigma(\g,x):=\inf\limits_{\lambda\in \hat\Sigma(\g,x)}\lambda, resp.\right).$$

We call $\sigma(\g,x)$ the $\sigma$-eigenvalue to $x$ for $\g$.
\dfd
\rmk
By the definition, it is clear that $\Sigma(\g,x)\subset\hat\Sigma(\g,x),\sigma(\g,x)\geq\hat\sigma(\g,x)$ if $\Sigma(\g,x)\neq\emptyset$.
\rmkd

Now we give another equivalent description of $\lambda$-flow satisfying (\ref{sigma2,2}) in Definition \ref{sigma2}.
\lm\label{positive flow}
Given any finite tree $\g=(V,E)$ with boundary $\dg$ containing $x$. Then $f_\lambda\in \rv$ is a $\lambda$-flow to $x$ with $\lambda>0$ satisfying (\ref{sigma2,2}) in Definition \ref{sigma2} if and only if $f_\lambda\in \rv$ is a $\lambda$-flow to $x$ with $\lambda>0$ satisfying $f_\lambda(z)>0$ for any $z\in \dg-\{x\}$.

\lmd
\pf
For any $z\in \dg-\{x\}$, there is a path $z_0=z\sim z_1\sim\cdots\sim z_i=x$ for some $i\in\N$. If $f_\lambda\in \rv$ is a $\lambda$-flow to $x$ satisfying Definition \ref{sigma2} with $\lambda\geq 0$, since $d(z,x)=i>i-1=d(z_1,x),$ we have $f_\lambda(z)-f_\lambda(z_1)=\lambda f_\lambda(z)>0.$ Thus we get $f_\lambda(z)>0$ as $\lambda>0$.

If $f_\lambda\in \rv$ is a $\lambda$-flow to $x$ with $\lambda>0$ satisfying $f_\lambda(z)>0$ for any $z\in \dg-\{x\}$, one easily checks that $f_\lambda$ satisfies (\ref{sigma2,2}) in Definition \ref{sigma2} by Lemma \ref{edge flow}.

%The case of $|V|=2$ is trivial.
%Let $x_1\in V-\dg,x_j\in V$ with $d(x,x_1)=d(x_1,x_j)=1$ and $2\leq j\leq m$ for $m=deg(x_1)$, namely $x,x_j$ are exactly all neighbours of $x_1$ for $2\leq j\leq m$. Denote by $\g_j=(V_j,E_j)$ the subtree $\h_k(x_1,x_j)$ with boundary $\dg_k$, where $\h_k$ is the branch from $(x_1,x_j)$. Note that $g_\lambda^k:=f_\lambda|_{V_k}$ is a $\lambda$-flow to $x_1$ on $\g_j$ with $\lambda>0$ satisfying $g_\lambda^k(z)>0$ for any $z\in \dg_k-\{x_1\}$. By induction, $g_\lambda^k$ satisfies (\ref{sigma2,2}) in Definition \ref{sigma2}, so we have $f_\lambda(x_j)-f_\lambda(x_1)=g^k_\lambda(x_j)-g^k_\lambda(x_1)=\nabla_{(x_j,x_1)} g^k_\lambda>0$. On the other hand, $\Delta f_\lambda(x_1)=0$, hence we obtain $f_\lambda(x_1)-f_\lambda(x)=\sum\limits_{j=1}^i\nabla_{(x_j,x_1)} f_\lambda>0$.

\pfd

\lm\label{po-steklov}
For a finite tree $\g=(V,E)$ with boundary $\dg$ containing $x$, there are at most $(|\dg|-1)$ many $\lambda$ satisfying (\ref{sigma2,1}) in Definition \ref{sigma2}.

As a consequence, we have $\sigma(\g,x)=\min\limits_{\lambda\in \Sigma(\g,x)}\lambda\geq\hat\sigma(\g,x)=\min\limits_{\lambda\in \hat\Sigma(\g,x)}\lambda\geq \lambda_2((\g)^2_x)>0$, if $\Sigma(\g,x)\neq \emptyset.$

\lmd
\pf
Assume $\lambda$ and $f_\lambda$ satisfy (\ref{sigma2,1}) in Definition \ref{sigma2}. Then we can define a function $g$ on the graph $\hat\g:=(\g)^2_x=(\hat V,\hat E)$ with boundary $\delta\hat{\mathcal{\g}}$, $\hat V=V\bigsqcup V'/x\backsim x'$, and $\hat E=E\bigsqcup E'$ as follows.
Let $g:\hat V\longrightarrow \R,\ z\longmapsto f_\lambda(z)$, if $z\in V$, $z'\longmapsto -f_\lambda(z)$ if $z'\in V'$.
Since $f_\lambda(x)=-f_\lambda(x)=0,$ the function $g$ is well-defined.

One easily verifies that $\Delta g(x)=0$ by Lemma \ref{edge flow}. Hence $g$ satisfies
\begin{equation*}%\label{Steklov 2}
\left\{
\begin{array}{ll}
\Delta g(z)=0, &z\in \hat\Omega:=\hat V-\delta\hat{\mathcal{\g}},\\
\dfrac{\partial g}{\partial n}(z)=\lambda g(z), &z\in \delta\hat{\mathcal{\g}}.
\end{array}
\right.
\end{equation*}
Note that $f_\lambda=0$ if $\lambda=0$, which is impossible by the definition of $\lambda$-flow. Thus we complete the proof.
\pfd

We introduce the key lemma on the existence and the uniqueness (up to scaling) of some $\lambda$-flow, which is continuous in $\lambda$ for sufficiently small $\lambda.$
\lm \label{exixt posi}
For any finite tree $\g=(V,E)$ with boundary $\dg$ containing $x$, assume $x_1,x_j\in V$ such that $x\sim x_1\sim x_j$ with $2\leq j\leq m:=\deg x_1$. Choose any vertex $w\in \dg-\{x\}$ and fix it. Then there exists the unique $\lambda$-flow to $x$, $f_\lambda\in \rv$ with $0\leq\lambda<\sigma_1$ such that
\begin{enumerate}[(1)]
\item (uniform sign and normalization) $f_\lambda(z)>0$ for all $z\in \dg-\{x\}$ with $f_\lambda(w)=1$,\label{p1}
%\item $\nabla_{(y,z)} f_\lambda=f_\lambda(y)-f_\lambda(z)>0$ if $\lambda>0$ and $\nabla_{(y,z)} f_\lambda=f_\lambda(y)-f_\lambda(z)=0$ if $\lambda=0$, where $(y,z)\in E$ with $d(x,y)>d(x,z)$.\label{p1}
\item (continuity) $f_\lambda$ is continuous in $\lambda$ for $0\leq\lambda<\sigma_1$, \label{p2}
\item (monotonicity) $0<\min\{0\leq\lambda<\sigma_1|f_\lambda(x)=0\}=\sigma(\g,x)<\sigma_1$, \label{p3}%\sigma(\g,x)>0.
%$f_\lambda(x)=0$ if and only if $\lambda=\sigma(\g,x)$?$\{0\leq\lambda<\sigma_1|f_\lambda(x)=0\}=\{0<\lambda<\sigma_1|f_\lambda(x)=0\}\neq \emptyset$, i.e.
\end{enumerate}
where $\sigma_1:=\min\limits_{2\leq j\leq m}\sigma(\g_j,x_1)$, $\g_j:=\h_j(x_j,x_1)=(V_j,E_j)$ with boundary $\dg_j$ and $\h_j$ is the branch from $(x_j,x_1)$.

Moreover, $f_\lambda$ is continuous in $\lambda$ as $0\leq\lambda\leq\sigma_1$ if there exists $j_0$ with $2\leq j_0\leq m$ such that $\sigma(\g_j,x_1)>\sigma(\g_{j_0},x_1)$, for any $2\leq j\neq j_0\leq m$.
\lmd

\pf
Denote by $\om_j:=V_j-\dg_j$ and $\om:=V-\dg$. Note that $x_1\in \dg_j\cap\om$ and that (\ref{p1}) holds is equivalent to the following holds that $\nabla_{(y,z)} f_\lambda=f_\lambda(y)-f_\lambda(z)>0$ if $\lambda>0$, where $(y,z)\in E$ with $d(x,y)>d(x,z)$ by Lemma \ref{positive flow}.

It is clear that such $f_\lambda$ is a constant for $\lambda=0$ and then $f_0=f_0(w)=1.$ So that $$\{0\leq\lambda<\sigma_1|f_\lambda(x)=0\}=\{0<\lambda<\sigma_1|f_\lambda(x)=0\}.$$

We argue by induction on $|V|$. The case of $|V|=2$ is trivial.

Since $x\notin V_j$, $|V_j|<|V|$. Suppose that $x_1,x_{j,l}$ are all neighbours of $x_j$ for $2\leq l\leq m_j$, where $m_j=\deg x_j$.  By induction, we have a $\lambda$-flow $f_{\lambda,j}$ to $x_1$ on $\g_j$ for $0\leq \lambda<\sigma_{1,j}$ such that
\begin{enumerate}[(i)]
\item  $f_{\lambda,j}(w_j)=1$ for fixed $w_j\in \dg_j-\{x_1\}$ and $f_{\lambda,j}(z)>0$ for all $z\in \dg_j-\{x_1\}$,\label{p2,1}
%\item $\nabla_{(y,z)} f_{\lambda,j}=f_{\lambda,j}(y)-f_{\lambda,j}(z)>0$ if $\lambda>0$ and $\nabla_{(y,z)} f_{\lambda,j}=f_{\lambda,j}(y)-f_{\lambda,j}(z)=0$ if $\lambda=0$, where $(y,z)\in E_j$ with $d(x,y)>d(x,z)$.
\item $f_{\lambda,j}$ is continuous in $\lambda$ for $0\leq\lambda<\sigma_{1,j}$,\label{p2,2}
\item $\{0\leq\lambda<\sigma_{1,j}|f_{\lambda,j}(x_1)=0\}=\{0<\lambda<\sigma_{1,j}|f_{\lambda,j}(x_1)=0\}\neq \emptyset$, i.e $0<\min\{0\leq\lambda<\sigma_{1,j}|f_{\lambda,j}(x_1)=0\}=\sigma(\g_j,x_1)<\sigma_{1,j},$\label{p2,3}
%$f_\lambda(x)=0$ if and only if $\lambda=\sigma(\g,x)$?
\end{enumerate}
where $\sigma_{1,j}:=\min\limits_{2\leq l\leq m_j}\sigma(\g_{j,l},x_j)$, $\g_{j,l}:=\h_{j,l}(x_{j,l},x_j)=(V_{j,l},E_{j,l})$ with boundary $\dg_{j,l}$ and $\h_{j,l}$ is the branch from $(x_{j,l},x_j)$.

First we will show $f_{\lambda,j}(x_1)>0$ for $0\leq \lambda<\sigma_1,2\leq j\leq m.$ Otherwise, $f_{\lambda_0,j_0}(x_1)\leq 0$ with some $0\leq \lambda_0<\sigma_1$ and $2\leq j_0\leq m$ and $f_{0,j_0}(x_1)=f_{0,j_0}(z)=f_{0,j_0}(w_{j_0})=1>0$ for all $z\in \dg_{j_0}-\{x_1\}$. By (\ref{p2,2}) in the above, there exists some $0<\hat\lambda_0\leq\lambda_0<\sigma_1$ such that $f_{\hat\lambda_0,j_0}(x_1)=0$. This is impossible by the definition of $\sigma_1$.

Note that $\sigma_1\leq \sigma(\g_j,x_1)<\sigma_{1,j}$. Thanks to the finiteness of $\g$, we may assume $\sigma_1=\sigma(\g_2,x_1)$. Then we construct $f_\lambda$ for $0\leq \lambda<\sigma_1$ as follows. Let
\begin{align}\label{consf}
f_\lambda(z):=c\dfrac{f_{\lambda,2}(x_1)}{f_{\lambda,j}(x_1)}f_{\lambda,j}(z)
\end{align}
for $z\in V_j$ and some $c>0$,
\begin{align*}
f_\lambda(x):=f_\lambda(x_1)-\sum\limits_{j=2}^m\nabla_{(x_j,x_1)}f_{\lambda}=cf_{\lambda,2}(x_1)-c\dfrac{f_{\lambda,2}(x_1)}{f_{\lambda,j}(x_1)}\sum\limits_{j=2}^m\nabla_{(x_j,x_1)}f_{\lambda,j},
\end{align*}
such that $f_\lambda(w)=1$.

One easily checks that $f_\lambda$ is well-defined and indeed it is a $\lambda$-flow to $x$ on $\g$ satisfying (\ref{p1}),(\ref{p2}) in Lemma \ref{exixt posi}. We need to show (\ref{p3}) in Lemma \ref{exixt posi}.

\begin{claim}\label{negative}
There exists some $\lambda_{0,1}\in (0,\sigma_1)$ such that $f_{\lambda_{0,1}}(x)<0$.
\end{claim}
\pf[Proof of Claim \ref{negative}]
%Let $\hat v\in \dg_2-\{x_1\}$ and $\hat v=v_0\sim v_1\sim\cdots\sim v_L=x_1$ be a path connecting $x_1$ and $\hat v$ in $V_j$, where $L=d(x_1,\hat v)$.

Recall that $\sigma_1=\sigma(\g_2,x_1)$ and for $\lambda\in [0,\sigma_1)$, $f_\lambda(x)=f_\lambda(x_1)-\sum\limits_{j=2}^m\nabla_{(x_j,x_1)}f_\lambda$, $f_\lambda(x_1)=cf_{\lambda,2}(x_1)$ and $\nabla_{(x_j,x_1)}f_\lambda>0$ by (\ref{p1}) in Lemma \ref{exixt posi}.

Then we have
$f_\lambda(x)< c(f_{\lambda,2}(x_1)-\nabla_{(x_2,x_1)}f_{\lambda,2})=:g_{\lambda,2}$.
By induction $g_{\lambda,2}$ is continuous in $\lambda$ for $0\leq \lambda<\sigma_{1,2},$ and hence continuous in $\lambda=\sigma(\g_2,x_1)=\sigma_1<\sigma_{1,2}$. Thus we obtain
$g_{\sigma_{1},2}=-c\nabla_{(x_2,x_1)}f_{\sigma(\g_2,x_1),2}<0$. The above discussion yields that $f_{\lambda_{0,1}}(x)<0$ for some $\lambda_{0,1}\in (0,\sigma_1)$.

\pfd
%\sum\limits_{n=1}^L\nabla_{(v_{n-1},v_{n})}f_{\lambda,t}
Note that $f_0=1,$ and then $f_0(x)>0$. By Claim \ref{negative} and (\ref{p2}), there exists some $\lambda_{0,2}$ such that $f_{\lambda_{0,2}}(x)=0$ with $0<\lambda_{0,2}\leq \lambda_{0,1}<\sigma_1$.
This gives $0<\sigma(\g,x)\leq\min\{0\leq\lambda<\sigma_1|f_\lambda(x)=0\}<\sigma_1$.

Now we prove the uniqueness. It suffices to show for any $g_\lambda$ satisfying conditions in Lemma \ref{exixt posi}, then $g_\lambda=f_\lambda.$ We also argue by induction on $|V|$, the case of $|V|=2$ is also trivial. By induction on $f_{\lambda,j},g_{\lambda,j}$ on $\g_j$ with $f_{\lambda,j}(w_j)=g_{\lambda,j}(w_j)=1$, where $2\leq j\leq m$, we have
\begin{align}
&f_{\lambda,j}=g_{\lambda,j},\label{fg1}
\\&f_\lambda|_{V_j}=c_1\dfrac{f_{\lambda,2}(x_1)}{f_{\lambda,j}(x_1)}f_{\lambda,j},\label{fg2}
\\&g_\lambda|_{V_j}=c_2\dfrac{g_{\lambda,2}(x_1)}{g_{\lambda,j}(x_1)}g_{\lambda,j}.\label{fg3}
\end{align}

Suppose $w\in V_i$ for some $i$ with $2\leq i\leq m$. On the other hand, by induction on $f_\lambda|_{V_i},g_\lambda|_{V_i}$ on $\g_i$ with $f_\lambda|_{V_i}(w)=g_\lambda|_{V_i}(w)=1$, we have $f_\lambda|_{V_i}=g_\lambda|_{V_i}$. Hence by (\ref{fg1}), (\ref{fg2}), (\ref{fg3}) for $j=i$ , we obtain $c_1=c_2$.
So we get $f_\lambda|_{V-\{x\}}=g_\lambda|_{V-\{x\}}$. Note that $$f_{\lambda}(x)=f_\lambda(x_1)-\sum\limits_{j=2}^m\nabla_{(x_j,x_1)}f_\lambda, g_{\lambda}(x)=g_\lambda(x_1)-\sum\limits_{j=2}^m\nabla_{(x_j,x_1)}g_\lambda.$$ It is clear that $f_{\lambda}(x)=g_{\lambda}(x)$.

This verifies (\ref{p3}) in Lemma \ref{exixt posi}.

Note that $\sigma_1=\sigma(\g_{j_0},x_1)$ if there exists $j_0$ with $2\leq j_0\leq m$ such that $\sigma(\g_j,x_1)>\sigma(\g_{j_0},x_1)$, for any $2\leq j\neq j_0\leq m$. Similarly, we can construct $f_\lambda$ by replacing $f_{\lambda,2}$ with $f_{\lambda,j_0}$ and we claim that $f_{\sigma_1,j}(x_1)>0$ for any $2\leq j\neq j_0\leq m$. Otherwise, $f_{\sigma_1,j}(x_1)\leq 0,f_{0,j}(x_1)=1>0$, then there exists some $\hat\sigma\in [0,\sigma_1]$ such that $f_{\hat\sigma,j}(x_1)=0$, by the continuity for $\lambda\in [0,\sigma(\g_j,x_1))$. So $\sigma(\g_j,x_1)\leq\hat\sigma\leq\sigma_1$, which contradicts that $\sigma(\g_j,x_1)>\sigma(\g_{j_0},x_1)=\sigma_1$. Hence we can construct $f_\lambda$ by replacing $f_{\lambda,2}$ with $f_{\lambda,j_0}$ in (\ref{consf}).

Finally, we conclude the proof.

%It is easy to verify that there exists $w\in \dg$ satisfying $d(x,w)=\max\limits_{v\in V}d(x,v)$. We may assume $d(x,w)=L$ and $x_1\in V$ with $d(x,x_1)=1,d(x_1,w)=L-1$.

%We argue by induction on $L$. The case of $L=1$ is trivial. For the case of $L=k+1\geq 2$, assume $deg(x_1)=m$ and $x_2,x_3,\cdots,x_m\in V$ with $d(x_1,x_2)=d(x_1,x_3)=\cdots=d(x_1,x_m)=1$. Denote by $\h_j$ the branch from $(x_1,x_j)$ and $\h_2$ contains $w$, where $2\leq j\leq m$.

%Note that $\max\limits_{w\in V_j}d(x_1,w)=L-1=k$, where $V_j$ is the vertices set of $\h_j(x_1,x_j)$. So by the inductive assumption, we have a $\lambda_j$-flow to $x_1$, $f_{\lambda_j}$ on the subtree $\h_j(x_1,x_j)$ with the relative boundary $\dg(\h_j(x_1,x_j))$ in $\g$, where $\lambda_j>0$. Moreover, for any edge $(u_j,v_j)$ in $\h_j(x_1,x_j)$ with $d(x_1,u_j)>d(x_1,v_j)$, $\nabla_{(u_j,v_j)} f_{\lambda_j}=f_{\lambda_j}(u_j)-f_{\lambda_j}(v_j)>0$, $f_{\lambda_j}(x_1)=0$ and $f_{\lambda_j}$ is positive on $V_j$. By Definition \ref{lambda-flow}, it is not hard to get $\nabla_{(u_j,v_j)} f_{\lambda_j}=\sum f_{\lambda_j}(w)\lambda_j$, where $w$ runs over the relative boundary vertics of the branch from $(u_j,v_j)$ in the summation.

%Let $x_{j,1}\in \dg(\h_j(x_1,x_j))$, we may assume that $x_1\sim y_{j,1}\sim y_{j,2}\cdots\sim x_{j_1}$

\pfd

As a consequence, we have the following.

\co\label{mono1}
Let $\h$ with boundary $\delta\h$ be a subtree of $\g=(V,E)$ with boundary $\dg$ containing $x.$ Assume that $y\in \delta\h\cap (V-\dg)$ and $\h$ does not contain $x$ such that $\h-y$ and $x$ are on the different sides of $y$. Then
$$\sigma(\h,y)>\sigma(\g,x).$$
\cod

\pf
 This follows from Lemma \ref{exixt posi}.
\pfd

\co\label{mono}
Let $f_\lambda(x)$ be as in Lemma \ref{exixt posi} and $f_\gamma(x)<0$ for some $0\leq \gamma<\sigma_1.$ Then $\gamma>\sigma(\g,x)$.

\cod
\pf
Note that $f_0=f_0(x)=1>0$ and $f_\gamma(x)<0$. Since $f_\lambda$ is continuous in $\lambda$ for $0\leq \lambda<\sigma_1$, one can deduce that
$f_{\hat\gamma}(x)=0$ for some $0<\hat\gamma<\gamma$. Thus $0<\sigma(\g,x)\leq \hat\gamma<\gamma.$
\pfd

Next, we propose a very useful criterion, which is used to compare $\lambda$ associated to some $\lambda$-flow on a finite tree $\g$ with $\sigma(\g,x)$, where $x$ is a boundary vertex of $\g$.

\lm\label{subtree}
Given any finite tree $\g=(V,E)$ with boundary $\dg$ containing $x$. Assume $f_\lambda$ is a $\lambda$-flow to $x$ on $\g$, then the following holds.

If $f_\lambda(z)=0$ for some $z\in \dg-\{x\}$ or $f_\lambda(z_1)f_\lambda(z_2)<0$ for some $z_1,z_2\in \dg-\{x\}$, then there exists a proper subtree $\g_1=(V_1,E_1)=\h(u,v)$ with boundary $\dg_1$ such that $x\notin V_1$, $\h$ and $x$ are on different sides of $v$, $f_\lambda(v)\leq0$, $f_\lambda(s)>0$ or $f_\lambda(v)\geq0$, $f_\lambda(s)<0$ for any $s\in \dg_1-\{v\}$, where $\h$ is the branch from $(u,v)$ in $\g$.
In particular, $\lambda>\sigma(\g,x)$.

\lmd

\pf
Recall that $f_\lambda|_{\dg-\{x\}}\neq 0$, we may assume that $f_\lambda(w)\neq 0$ for some $w\in \dg-\{x\}.$

Case 1. $f_\lambda(z)=0$ for some $z\in\dg-\{x\}.$ We may assume that there exist $y$ and its neighbors $y_1,y_2,\cdots,y_k$ with $k=\deg(y)$ such that $z\in \h_1(y_1,y),w\in \h_2(y_2,y),x\in \h_k(y_k,y)$, where $\h_j$ is the branch from $(y_j,y)$ with $1\leq j\leq k$.

If $f_\lambda|_{\delta\h_1(y_1,y)-\{y\}}\neq 0$, we get the proof by induction for $\dh_1(y_1,y)$. Otherwise, $f_\lambda|_{\delta\h_1(y_1,y)-\{y\}}=0$. It is clear that $f_\lambda|_{\h_1(y_1,y)}=0$ by the definition of $f_\lambda$. In particular, we have $f_\lambda(y)=0$.

If $f_\lambda(z^\prime)f_\lambda(w)\leq 0$ for some $z^\prime\in \dh_2(y_2,y)-\{y\}$, we get the proof by induction for $\dh_2(y_2,y)$. Otherwise, $f_\lambda(s)f_\lambda(w)> 0$ for any $s\in\dh_2(y_2,y)-\{y\}$. Recall that $f_\lambda(y)=0$ and $\h_2(y_2,y)-\{y\}$ and $x$ are on different sides of $y$, we get the proof.

Case 2. $f_\lambda(z_1)f_\lambda(z_2)<0$ for some $z_1,z_2\in \dg-\{x\}$. The argument is same to that for Case 1.
\pfd

\co\label{mono3}
Given a finite tree $\g=(V,E)$ with boundary $\dg$ containing $x$, then $\sigma(\g,x)=\hat\sigma(\g,x).$
\cod
\pf
%Let $f_\lambda\in \rv$ be a $\lambda$-flow to $x$ on $\g$.
%If there exist $v_1,v_2\in \dg-\{x\},v_3\in V-\dg$ such that $g_2(v_1)g_2(v_2)\leq 0$, one can suppose $d(v_3,v_1)=1,g_2(v_1)\leq0,g_2(v_2)>0$. We may assume  $f_\lambda(v_3)\leq0$, otherwise it follows that $\lambda\geq1$ from $f_\lambda(v_3)=f_\lambda(v_1)(1-\lambda(\g))$, which is trivial. Using Lemma \ref{subtree} we can find a proper subtree $\h=(V_1,E_1)$ in $\g$ with boundary $\delta\h$ containing $v_2$, $y\in \delta\h\cap V-\dg$ and $v_1\notin V_1$ such that $f_\lambda(y)<0$, $f_\lambda(z)>0$ for $z\in \delta\h-\{y\}$.
This is a direct consequence of Lemma \ref{mono1}, Lemma \ref{subtree}.

\pfd

%\lambda-flow
\section{Proof of Theorem \ref{main1}}\label{sec:proof}

Before giving the proof of Theorem \ref{main1}, we will give the relation between $\lambda_2(\g)$ and $\sigma$-eigenvalues for the branches in $\g$.

\tm\label{sigmasteklov}
 Given any finite tree $\g=(V,E)$ with boundary $\dg$, assume $x,x_j\in V$ such that $d(x,x_j)=1$ with $1\leq j\leq m:=\deg(x)$. Denote by $\g_j:=\h_j(x_j,x)=(V_j,E_j)$ with boundary $\dg_j$, where $\h_j$ is the branch from $(x_j,x)$ in $\g$. Assume that $\sigma(\g_1,x)\leq \sigma(\g_2,x)\leq \cdots\leq\sigma(\g_m,x)$, then the following statements holds.

 If $\sigma(\g_1,x)=\sigma(\g_2,x)$, then $\lambda_2(\g)=\sigma(\g_1,x)$. Moreover, any $\lambda_2(\g)$-eigenfunction must be zero on $x$.

 If $\sigma(\g_1,x)<\sigma(\g_2,x)$, then $\sigma(\g_1,x)<\lambda_2(\g)<\sigma(\g_2,x)$.
\tmd
\pf
Let $f_{\lambda_2(\g)}$ be any Steklov eigenvalue with eigenvalue $\lambda_2(\g)$. Since $\sum_{y\in \dg}f_{\lambda_2(\g)}(y)=0$, we may assume $f_{\lambda_2(\g)}(u)f_{\lambda_2(\g)}(x)\leq 0$, where $u\in V_j$ for some $1\leq j\leq m$. Using Lemma \ref{subtree} to the branch $\g_j$, we have 
\begin{align}\label{zero}
\lambda_2(\g)>\sigma(\g_j,x)\geq\sigma(\g_1,x)
\end{align}
if $f_{\lambda_2(\g)}(x)\neq 0$. On the other hand, it is clear that $\lambda_2(\g)\geq\sigma(\g_j,x)\geq \sigma(\g_1,x)$ via Lemma \ref{mono3} if $f_{\lambda_2(\g)}(x)=0$. So we get $\lambda_2(\g)\geq\sigma(\g_1,x).$

Now one can construct a function as follows. Set
\begin{equation*}%\label{Steklov 2}
\left\{
\begin{array}{ll}
 g(z)=af_{\sigma(\g_1,x)}(z), & z\in V_1,\\
g(z)=bf_{\sigma(\g_2,x)}(z),& z\in V_2,\\
g(z)=0,& elsewhere,
\end{array}
\right.
\end{equation*}
where $f_{\sigma(\g_1,x)}$ ($f_{\sigma(\g_2,x)}$, resp.) is the $\sigma(\g_1,x)$ ($\sigma(\g_2,x)$, resp.) flow to $x$ on $\g_1$ ($\g_2$, resp.) as in Lemma \ref{exixt posi} and $a=\sum_{y\in \delta\g_2}f_{\sigma(\g_2,x)}(y),b=-\sum_{y\in \delta\g_1}f_{\sigma(\g_1,x)}(y)$.

Note that $\sum_{y\in\dg}g(y)=0$ and $R(g)=\sigma(\g_1,x)$ if $\sigma(\g_1,x)=\sigma(\g_2,x)$ and $\sigma(\g_1,x)<R(f)<\sigma(\g_2,x)$ if $\sigma(\g_1,x)<\sigma(\g_2,x)$. This yields $\lambda_2(\g)=\sigma(\g_1,x)$ if $\sigma(\g_1,x)=\sigma(\g_2,x)$ and $\sigma(\g_1,x)\leq\lambda_2(\g)<\sigma(\g_2,x)$ if $\sigma(\g_1,x)<\sigma(\g_2,x)$. Recall the inequality (\ref{zero}), we obtain $f_{\lambda_2(\g)}(x)=0$ if $\sigma(\g_1,x)=\sigma(\g_2,x).$

If $\sigma(\g_1,x)<\sigma(\g_2,x)$ and  $\lambda_2(\g)=\sigma(\g_1,x)$, we obtain $f_{\lambda_2(\g)}(x)=0$ and $f_{\lambda_2(\g)}(y)\neq 0$ for some $y\in V_{j_0}$ with $2\leq j_0\leq m$. Using Corollary \ref{mono3} to $\g_{j_0}$, we get $\lambda_2(\g)\geq\sigma(\g_{j_0},x)\geq\sigma(\g_2,x)>\sigma(\g_1,x)=\lambda_2(\g)$, which is a contradiction.
\pfd

\co\label{partition}
Given any finite tree $\g=(V,E)$ with boundary $\dg$ and $\hat{\mathcal{\g}}:=(\g)_x^2=(\hat V,\hat E)$ with boundary $\delta\hat\g$, then the following hold.
\begin{enumerate}[1)]
\item $\lambda_2(\hat{\mathcal{\g}})=\sigma(\g,x)<\lambda_2(\g)$ if $x\in \dg$.
\item $\lambda_2(\hat{\mathcal{\g}})\leq\lambda_2(\g)$ if $x\in V-\dg$.
\end{enumerate}
\cod
\pf
These are direct consequences of Corollary \ref{mono3} and Theorem \ref{sigmasteklov}.

\pfd

%{As a consequence, $\lambda_2(\g)\geq \min\limits_{1\leq j\leq \deg(z)}\sigma(\f_j(z_j,z),z)$.}

%\pf[Proof of Claim \ref{parti1}]

Now we are ready to prove Theorem \ref{main1}.

\pf[Proof of Theorem \ref{main1}]
We will first give the proof of monotonicity, and then show the rigidity based on the proof of monotonicity.

It suffices to show that for the case that two finite trees $\g=(V,E)\subset\g^\prime=(V^\prime,E^\prime)$ such that $V^\prime=V\cup\{y\},E^\prime=E\cup (x,y)\cup(y,x)$ for $x\in V,y\in V^\prime$.

Denote by $\dg,\dg^\prime$ the boundary of $\g,\g^\prime$ respectively, and by $\om:=V-\dg,\om^\prime:=V^\prime-\dg^\prime$ the set of interior vertices respectively. Let $h_2\in \rv$ be a Steklov eigenfunction associated with $\lambda_2(\g)$. Assume $h_2(x_0)>0$ for some $x_0\in \dg$.

\textbf{The proof of the monotonicity.}

%\red{Step 1. Find an directed edge $(x_k,x_{k+1})$ such that $h_2|\g_k>0$ and we can assume $h_2|\g_k^\prime\leq 0$, where $x_k,x_{k+1}\in V$ and $\g_k$ is the branch from $(x_k,x_{k+1})$ in $\g$ and $\g_k^\prime$ is the branch from $(x_{k+1},x_k)$ in $\g$.}

Let $x_1\in V$ with $d(x_0,x_1)=1$. Note that $1=\sigma(\g_0(x_0,x_1),x_1)\geq\lambda_2((\g)^2_{x_1})$ by Corollary \ref{partition}, where $\g_0=(V_0,E_0)$ is the branch from $(x_0,x_1)$ in $\g$, i.e. $\{x_0\}$.

If $1=\sigma(\g_0(x_0,x_1),x_1)\leq\lambda_2(\g)$, then $\sigma(\g_0(x_0,x_1),x_1)=\lambda_2(\g)=1$ and this case is trivial.

If $\sigma(\g_0(x_0,x_1),x_1)>\lambda_2(\g)$ and $\sigma(\g_1^\prime(x_2,x_1),x_1)=\lambda_2((\g)^2_{x_1})$ for some $x_2\in V$ with $d(x_1,x_2)=1$, where $\g_1^\prime$ is the branch from $(x_2,x_1)$.
Since $\g$ is finite, using Corollary \ref{mono1}, Theorem \ref{sigmasteklov}, and Corollary \ref{partition}, we can find $x_j\in V$ such that $\sigma(\g_{j-1}(x_{j-1},x_{j}))>\lambda_2(\g),$ $\sigma(\g_j^\prime(x_{j+1},x_j),x_j)=\lambda_2((\g)^2_{x_j})$ for $1\leq j\leq k,$ and
\begin{align*}
\sigma(\g_k(x_k,x_{k+1}),x_{k+1})=\lambda_2((\g)^2_{x_{k+1}})\leq\lambda_2(\g),
\end{align*}
where $\g_j=(V_j,E_j)$ is the branch from $(x_j,x_{j+1})$ and $\g_j^\prime=(V_j^\prime,E_j^\prime)$ is the branch from $(x_{j+1},x_j)$, $k\geq 2$.

In summary, we have
\begin{align}
&\sigma(\g_{k-1}(x_{k-1},x_{k}))>\lambda_2(\g),\label{k-1geq}
\\&\sigma(\g_k^\prime(x_{k+1},x_k),x_k)\leq\lambda_2(\g),\label{kprimeleq}
\\&\sigma(\g_k(x_k,x_{k+1}),x_{k+1})\leq\lambda_2(\g).\label{kleq}
\end{align}

Let $x_{k,i},x_{k+1,t}\in V$ such that $d(x_{k,i},x_k)=d(x_{k+1,t},x_{k+1})=1$ with $1\leq i\leq \deg(x_k)-2,1\leq t\leq \deg(x_{k+1})-1.$ Denote by $\g_{k,i}$ the branch from $(x_{k,i},x_k)$ in $\g,$ and by $\g_{k+1,t}=(V_{k+1,t},E_{k+1,t})$ the branch from $(x_{k+1,t},x_{k+1})$ in $\g$.

Case 1. Subcase 1. $\sigma(\g_k(x_{k},x_{k+1}),x_{k+1})=\lambda_2(\g)$. By Theorem \ref{sigmasteklov}, we have
\begin{align}
\sigma(\g_{k+1,t}(x_{k+1,t},x_{k+1}),x_{k+1})\geq \lambda_2(\g).\label{k+1geq}
\end{align}

Subcase 2. $\sigma(\g_k^\prime(x_{k+1},x_{k}),x_{k})=\lambda_2(\g)$. By Theorem \ref{sigmasteklov}, we have
\begin{align}
 \sigma(\g_{k,i}(x_{k,i},x_{k}),x_{k})\geq \lambda_2(\g).\label{kgeq}
\end{align}

Note that Subcase 1 is dual to Subcase 2. So that it suffices to consider Subcase 1.

By Theorem \ref{sigmasteklov}, there is $x_{k+1,t_0}\in V$ with $1\leq t_0\leq\deg(x_{k+1})-1$, such that $\sigma(\g_{k+1,t_0}(x_{k+1,t_0},x_{k+1}),x_{k+1})=\lambda_2(\g).$

Subcase 2.1.   If $x\notin V_k\cup V_{k+1,t_0},$ recalling that $\sigma(\g_{k}(x_{k},x_{k+1}),x_{k})=\lambda_2(\g),$ similar to the proof of Theorem \ref{sigmasteklov}, we can construct $g\in\rv$ with $g(x_0)=1$ such that
\begin{equation}\label{g}
g(z)=\left\{
\begin{array}{ll}
\sum\limits_{s\in \delta\g_{k+1,t_0}(x_{k+1,t_0},x_{k+1})-\{x_{k+1}\}}g_1(s)g_2(z), &z\in V_k,\\
-\sum\limits_{s\in \delta\g_{k}(x_k,x_{k+1})-\{x_{k+1}\}}g_2(s)g_1(z),&z\in V_{k+1,t_0},\\
0,&  otherwise,
\end{array}
\right.
\end{equation}
where $g_1,g_2$ are $\lambda_2(\g)$-flows to $x_{k+1}$ on $\g_{k+1,t_0}(x_{k+1,t_0},x_{k+1})$ and $\g_k(x_k,x_{k+1})$ respectively. One easily verifies that $g$ is a Steklov eigenfunction on $\g$ associated with eigenvalue $\lambda_2(\g)$. It is clear that $g(x)=0$ by $x\notin V_k\cup V_{k+1,t_0}.$ Then we extend $g$ to $g^\prime$ on $\g^\prime$ by setting $g^\prime(y)=0$. Note that $\sum\limits_{s\in \dg^\prime}g^\prime(s)=\sum\limits_{s\in \dg}g(s)=0$. Using the variational inequality (\ref{var k}) and Lemma \ref{scale}, we get
$$\lambda_2(\g^\prime)\leq R(g^\prime)=R(g)=\lambda_2(\g).$$

Subcase 2.2.   If $x\in V_k\cup V_{k+1,t_0},$ it suffices to show the proof for
\begin{align}\label{tglamdda}
\sigma(\g_{k+1,t}(x_{k+1,t},x_{k+1}),x_{k+1})>\lambda_2(\g)
\end{align}
with $1\leq t\neq t_0\leq\deg(x_{k+1})-1$, otherwise it can be reduced to Subcase 2.1. We may assume that $x\in V_{k+1,t_0}$, then there is a path $\alpha:=x_{k+1}=v_0\sim x_{k+1,t_0}=v_1\sim v_2\sim\cdots\sim v_n=x\sim y$ in $\g^\prime$ for some $n\geq 1$.%x=x_{k+1}

We will construct inductively a $\lambda$-flow to $y$ on $\g^\prime$ along the path $\alpha$, such that $f_\lambda$ is continuous in $\lambda$ as $0\leq \lambda\leq \lambda_2(\g)$ and $f_\lambda(x_0)=1$.

By (\ref{tglamdda}) and $\sigma(\g_k(x_{k},x_{k+1}),x_{k+1})=\lambda_2(\g)$, it follows from Lemma \ref{exixt posi} that there exists a $\lambda$-flow $f_{\lambda,1}$ to $v_1$ on $\h_1(v_0,v_1)$ such that $f_{\lambda,1}$ is continuous in $\lambda$ as $0\leq \lambda\leq \lambda_2(\g)$ and $f_{\lambda,1}(x_0)=1$, where $\h_1=(\hat V_{1,0},\hat E_{1,0})$ is the branch from $(v_0,v_1)$ in $\g$. By induction, there is a $\lambda$-flow $f_{\lambda,i}$ to $v_i$ on $\h_i(v_{i-1},v_i)$, such that $f_{\lambda,i}$ is continuous in $\lambda$ as $0\leq \lambda\leq \lambda_2(\g)$, where $\h_i=(\hat V_{i,i-1},\hat E_{i,i-1})$ is the branch from $(v_{i-1},v_i)$ in $\g$ for some $1\leq i\leq n$.

Assume that $v_{i,j}\sim v_i$ for $1\leq j\leq \deg(x_i)-2.$ Note that
\begin{align}\label{iglam}
\sigma(\h_{i,j}(v_{i,j},v_{i}),v_{i})> \sigma(\g_{k+1,t_0}(x_{k+1,t_0},x_{k+1}),x_{k+1})=\lambda_2(\g),
\end{align}
which follows from $\h_{i,j}(v_{i,j},v_{i})\subsetneq\g_{k+1,t_0}(x_{k+1,t_0},x_{k+1})$ and Lemma \ref{exixt posi}, where $\h_{i,j}=(\hat V_{i,j,i},\hat E_{i,j,i})$ is the branch from $(v_{i,j},v_{i})$ in $\g$. Then by Lemma \ref{exixt posi}, there is a $\lambda$-flow $f_{\lambda,i,j}$ to $v_i$ such that $f_{\lambda,i,j}$ is continuous in $\lambda$ as $0\leq \lambda\leq \lambda_2(\g)$ and $f_{\lambda_2(\g)}(v_i)>0$. We define $f_{\lambda,i+1}$ as follows.
\begin{equation}
f_{\lambda,i+1}(z)=\left\{
\begin{array}{ll}
f_{\lambda,i}(z), &z\in\hat V_{i,i-1},\\
\dfrac{f_{\lambda,i}(v_i)}{f_{\lambda,i,j}(v_i)}f_{\lambda,i,j}(z),&z\in \hat V_{i,j,i},\\
f_{\lambda,i}(v_i)-\nabla_{(v_{i-1},v_i)}f_{\lambda,i}-\sum\limits_{j=1}^{\deg(v_i)-2}\nabla_{(v_{i,j},v_i)}\dfrac{f_{\lambda,i}(v_i)}{f_{\lambda,i,j}(v_i)}f_{\lambda,i,j},&z=v_{i+1}. \\
\end{array}
\right.
\end{equation}
One easily checks that $f_{\lambda,i+1}$ is a $\lambda$-flow to $v_{i+1}$ on $\h_{i+1}(v_i,v_{i+1})$ such that $f_{\lambda,i+1}$ is continuous in $\lambda$ as $0\leq \lambda\leq \lambda_2(\g)$ and $f_{\lambda,i+1}(x_0)=1$.

Thus we construct a $\lambda$-flow to $y$ on $\g^\prime$ along the path $\alpha$, such that $f_\lambda$ is continuous in $\lambda$ as $0\leq \lambda\leq \lambda_2(\g)$ and $f_{\lambda}(x_0)=1$.

By the construction of $f_\lambda$, it is clear that $f_\lambda(v_0)>0$ as $0\leq \lambda<\lambda_2(\g)$, $f_{\lambda_2(\g)}(v_0)=0.$ Then by Lemma \ref{edge flow} or Lemma \ref{positive flow}, $f_{\lambda_2(\g)}(v_1)<0.$  Recalling (\ref{iglam}) and by the construction of $f_\lambda$, we have
\begin{align}\label{xleq0}
f_{\lambda_2(\g)}|_{ V_{k+1,t_0}}<0.
\end{align}

Note that since $f_{\lambda_2(\g)}(v_0)=0,$ $f_{\lambda_2(\g)}|_{V-V_k-\hat V_{1,0}}=0$. It is not hard to see that $f_{\lambda_2(\g)}|_{V}=g$ by the construction of $f_\lambda$ and Lemma \ref{exixt posi}, where $g$ is defined in (\ref{g}) and it is a Steklov eigenfunction associated with $\lambda_2(\g)$. So that we have
\begin{align}
&\sum\limits_{z\in\dg}f_{\lambda_2(\g)}(z)=0, \label{sum0}
\\&f_{\lambda_2(\g)}(y)=f_{\lambda_2(\g)}(x),\ if\ x\in\om,\label{xom}
\\&f_{\lambda_2(\g)}(y)=f_{\lambda_2(\g)}(x)-\nabla_{(x,y)}f_{\lambda_2(\g)},\ if\ x\in\dg.\label{xdo}
\end{align}
Hence we have
\begin{align}
&\sum\limits_{z\in\dg^\prime}f_{\lambda_2(\g)}(z)=f_{\lambda_2(\g)}(x),\ if\ x\in\om,\label{xom1}
\\&\sum\limits_{z\in\dg^\prime}f_{\lambda_2(\g)}(z)=-\nabla_{(x,y)}f_{\lambda_2(\g)}=\lambda_2(\g)f_{\lambda_2(\g)}(x),\ if\ x\in\dg.\label{xdo1}
\end{align}

Using (\ref{xom1}), (\ref{xdo1}) and (\ref{xleq0}), we have $\sum\limits_{z\in\dg^\prime}f_{\lambda_2(\g)}(z)<0$. Note that $\sum\limits_{z\in\dg^\prime}f_{0}(z)\geq f_0(x_0)>0$ and $f_\lambda$ is continuous in $\lambda$ for $0\leq \lambda\leq \lambda_2(\g).$ Thus there exists a $\hat\lambda$ such that $\sum\limits_{z\in\dg^\prime}f_{\hat\lambda}(z)=0$ for $0<\hat\lambda<\lambda_2(\g)$. Then $f_{\hat\lambda}$ is a Steklov eigenfunction on $\g^\prime$ associated with eigenvalue $\hat\lambda.$
By the inequality (\ref{var k}), we have
$$\lambda_2(\g^\prime)\leq R(f_{\hat\lambda})=\hat\lambda<\lambda_2(\g).$$

Case 2.
$\sigma(\g_{k-1}(x_{k-1},x_{k}))>\lambda_2(\g)$, $\sigma(\g_k^\prime(x_{k+1},x_k),x_k)<\lambda_2(\g)$, and $\sigma(\g_k(x_k,x_{k+1}),x_{k+1})<\lambda_2(\g).$
By Theorem \ref{sigmasteklov}, we have $$\sigma(\g_{k+1,t}(x_{k+1,t},x_{k+1}),x_{k+1})>\lambda_2(\g), \quad\sigma(\g_{k,i}(x_{k,i},x_k))>\lambda_2(\g),$$ for any $1\leq i\leq \deg(x_k)-2$ and $1\leq t\leq \deg(x_{k+1})-1$. We assume $x\in V_{k+1,t_0}\cup \{x_{k+1}\}$ and there is a path $\alpha:=v_0=x_{k+1}\sim x_{k+1,t_0}=v_1\sim\cdots\sim v_n=x\sim y$ or $x=x_{k+1}\sim y$. Similar to Case 1, we can construct a $\lambda$-flow $f_\lambda$ to $y$ on $\g^\prime$ with $f_\lambda(x_0)=1$ along the path $\alpha$, such that $f_{\lambda_2(\g)}$ is continuous in $\lambda$ for $0\leq\lambda\leq\lambda_2(\g)$.

This case is similar to Case 1. In fact, it is easier than Case 1, since $f_\lambda(v_i)\neq 0$ for any $\lambda\in [0,\lambda_2(\g)]$ and $v_i$. One easily checks that $f_{\lambda_2(\g)}|_{V}$ is a Steklov eigenfunction associated with the eigenvalue $\lambda_2(\g)$ and $f_{\lambda_2(\g)}(x)<0.$ Thus we have $\lambda_2(\g^\prime)<\lambda_2(\g)$.

\textbf{The proof of the rigidity.}

Now we will show the sufficient and necessary condition of the equality.

The ``if" part follows from Theorem \ref{sigmasteklov}. For the ``only if" part, i.e. $\lambda_2(\g^\prime)=\lambda_2(\g)$. By the above discuss, it only occurs in Case 1. We may assume that $\sigma(\g_{k}(x_{k},x_{k+1}),x_{k+1})=\sigma(\g_{k+1,t_0}(x_{k+1,t_0},x_{k+1}),x_{k+1})=\lambda_2(\g).$ Hence at least one of $$\sigma(\g^\prime_{k}(x_{k},x_{k+1}),x_{k+1})=\lambda_2(\g)$$
and
$$\sigma(\g^\prime_{k+1,t_0}(x_{k+1,t_0},x_{k+1}),x_{k+1})=\lambda_2(\g)$$
holds and the ``only if" part follows from Theorem \ref{sigmasteklov}.
\pfd

Now we have the following descriptions for the rigidity.
\co\label{critical}
Given any two finite trees $\g_1\subset\g_2$, suppose that $\lambda_2(\g_1)=\lambda_2(\g_2)$ and $\g_1\neq\g_2$. Then there exists a vertex $x$ in $\g_1$ such that any $\lambda_2(\g_1)$-eigenfunction, $\lambda_2(\g_2)$-eigenfunction on $\g_1,\g_2$ must be zero on $x$.
\cod
\pf
This is clear from Theorem \ref{sigmasteklov} and Theorem \ref{main1}.
\pfd

\co\label{strongrigidity}
Given any two finite trees $\g_1,\g_2$, suppose that there exist two adjacent vertices $x\sim y$ in $\g_1$ such that $\sigma(\f(y,x),x)<\lambda_2(\g_1)$ and
$\sigma(\h(x,y),y)<\lambda_2(\g_1)$,where $\f,\h$ are the branches from $(y,x),(x,y)$ in $\g_1$ respectively.
If $\g_1\subset \g_2$ or $\g_2\subset \g_1$, then $\lambda_2(\g_1)=\lambda_2(\g_2)$ if and only if $\g_1=\g_2$.
\cod

\pf
One can easily obtain that $\sigma(\f_i(x_i,x),x)>\lambda_2(\g_1),\sigma(\h_j(y_j,y),y)>\lambda_2(\g_1)$ for $2\leq i\leq \deg(x),2\leq j\leq \deg(y)$ by Theorem \ref{sigmasteklov}, where $x_i\sim y,y_j\sim x$ and $\f_i,\g_j$ are the branches from $(x_i,x),(y_j,y)$ in $\g_1$ respectively. If $\g_1\neq\g_2$, then there exist $z$, its neighbors $z_1,z_2$ in $\g_1$ such that
\begin{align}\label{strongrigidity1}
\sigma(\f^1(z_1,z),z)=\sigma(\f^2(z_2,z),z)=\lambda_2(\g_1),
\end{align}
where $\f^1,\f^2$ are the branches from $(z_1,z),(z_2,z)$ in $\g_1$ respectively. But (\ref{strongrigidity1}) is impossible for $z=x,y$. If $z\neq x,y$, then $z$ must in some branch from $\f_i,\g_j$ for $2\leq i\leq \deg(x),2\leq j\leq\deg(y)$, and (\ref{strongrigidity1}) is impossible by Corollary \ref{mono1}.

\pfd

The subsequent result implies some strict monotonicity of $\sigma(\g,x)$ for finite trees, provided that $x$ is a boundary vertex.

\co\label{smono}
Assume that $\g_1=(V_1,E_1)$ is a finite tree with boundary $\delta\g_1$ containing $x$, $\g_2=(V_2,E_2)$ is a proper subgraph of $\g_1$ (i.e. $V_2\subsetneq V_1$) such that $x\in V_2$. Then
\begin{align}\label{smo1}
\sigma(\g_2,x)>\sigma(\g_1,x).
\end{align}
\cod

\pf
It suffices to show (\ref{smo1}) holds for the case that $V_1=V_2\cup\{y\},E_1=E_2\cup\{(z,y),(y,z)\},x,z\in V_2$ and $y\notin V_2$. Consider three trees $\f_1:=(\g_2)^2_x\subset \f_2:=\g_1\bigsqcup\g_2/x\sim x\subset \f_3:==(\g_1)^2_x$, one can have $\sigma(\g_2,x)\geq \sigma(\g_1,x)$ and $\sigma(\g_2,x)=\lambda_2(\f_1)>\lambda_2(\f_2)\geq\sigma(\g_1,x)$ by Theorem \ref{sigmasteklov} , Theorem \ref{main1} and Subcase 2.2 of the proof of Theorem \ref{main1}.

\pfd

\section{Some applications}
%\cite{He,Hua}
In the previous paper \cite{HH20}, we proved that the bound for the first non-zero Steklov eigenvalue on any finite tree, in terms of the diameter of the tree.
\tm[Theorem 1.4, \cite{HH20}]\label{diames}
For any finite tree $\g=(V,E)$,
$$\lambda_2(\g)\leq\dfrac{2}{L},$$
where $L$ is the diameter of $\g$.
\tmd

Now we give a complete description for the case that attains the upper bound in Theorem \ref{diames}.
\tm\label{optimal diam}
Given a finite tree $\g=(V,E)$, assume there is a path $\h$ $x=x_0\sim x_1\sim x_2\cdots\sim x_L=y$ in $\g$, where $L$ is the diameter of $\g$. %For $x_k$ with $1\leq k\leq L-1$, there are three connected components if one removes edges $(x_{k-1},x_k),(x_k,x_{k+1})$.
 For $1\leq k\leq L-1,$ we remove the edges $(x_{k-1},x_k)$ and $(x_k,x_{k+1})$ from $\g,$ and denote by $\g_k=(V_k,E_k)$ the connected component containing $x_k,$ and by $\delta(\g_k,\g):=V_k\cap\dg$ the relative boundary. Set $n_k:=|\delta(\g_k,\g)|\in\N$.

Then $\lambda_2(\g)=\dfrac{2}{L}$ if and only if $L$ is odd, $n_k=0$ for $1\leq k\leq L-1$, i.e. it is just a path with length $L$ or $L$ is even, $n_k=0$ for $1\leq k\neq \dfrac{L}{2}\leq L-1$, and $\g=\h_1\cup\h_2$ such that $\lambda_2((\h_2)^2_{x_{\frac{L}{2}}})\geq \dfrac{2}{L}$, where the subgraph $\h_1=(V_1,E_1)$ is generated by all vertices $x_j$ with $0\leq j\leq L$ in $\g$ and the subgraph $\h_2=\g_{\frac{L}{2}}$.

\tmd

\pf
Note that $\g\supset\h$.

($\Rightarrow$)
If $L$ is odd, then $\sigma(\f_1(x_\frac{L-1}{2},x_\frac{L+1}{2}),x_\frac{L+1}{2})=\sigma(\f_2(x_\frac{L+1}{2},x_\frac{L-1}{2}),x_\frac{L-1}{2})=\dfrac{2}{L+1}<\dfrac{2}{L}=\lambda_2(\h)$, where $\f_1,\f_2$ are the branches from $(x_\frac{L-1}{2},x_\frac{L+1}{2}),(x_\frac{L+1}{2},x_\frac{L-1}{2})$ in $\h$ respectively. Thus $\g=\h$ by Corollary \ref{strongrigidity}.

If $L$ is even, we have $a_k=0$ for any $1\leq k\neq \dfrac{L}{2}\leq L-1$. Otherwise, note that $\sigma(\f_1(x_\frac{L-2}{2},x_\frac{L}{2}),x_\frac{L}{2})=\sigma(\f_2(x_\frac{L+2}{2},x_\frac{L}{2}),x_\frac{L}{2})=\dfrac{2}{L}$, where $\f_1,\f_2$ are the branches from $(x_\frac{L-2}{2},x_\frac{L}{2}),(x_\frac{L+2}{2},x_\frac{L}{2})$ in $\h$ respectively. We may assume that $a_k\neq 0$ for some $k\leq \frac{L-2}{2}$. Then by Corollary \ref{mono3} we have
$$\sigma(\f^1(x_\frac{L-2}{2},x_\frac{L}{2}),x_\frac{L}{2})<\sigma(\f_1(x_\frac{L-2}{2},x_\frac{L}{2}),x_\frac{L}{2})$$
and
$$\sigma(\f^2(x_\frac{L+2}{2},x_\frac{L}{2}),x_\frac{L}{2})\leq\sigma(\f_2(x_\frac{L+2}{2},x_\frac{L}{2}),x_\frac{L}{2}),$$
where $\f^1,\f^2$ are the branches from $(x_\frac{L-2}{2},x_\frac{L}{2}),(x_\frac{L+2}{2},x_\frac{L}{2})$ in $\g$ respectively. This is impossible by Theorem \ref{sigmasteklov} and we get the proof by Theorem \ref{main1}.

($\Leftarrow$)
This is straightforward by Theorem \ref{main1}.

\pfd

We denote by $T_{n+1}$ the $(n+1)$-regular tree for $n\in \N.$ We write $B_R(x)$ in a tree $\g$ for the ball centered at the vertex $x$ with radius $R$, with respect to the canonical combinatorial metric in $\g$.

Another interesting consequence is an explicit lower bound of the first Steklov eigenvalue for finite trees, with bounded vertex degree and bounded diameter.
Moreover, we give the characterization for the cases attaining the lower bound.

\tm\label{large tree}
Let $\g=(V,E)$ be a finite tree with diameter at most $L$ and vertex degree at most $D+1,$ i.e. $\deg(v)\leq D+1$ for any $v\in V$. Then the following hold.
\begin{enumerate}[(i)]
\item If $L=2R$ for some positive integer $R$, then
\begin{align}\label{bigtree1}
\lambda_2(\g)\geq \dfrac{1}{\sum\limits_{i=0}^{R-1}D^i}=\dfrac{D-1}{D^R-1}.
\end{align}
Moreover, the equality holds if and only if $\g$ contains a subgraph $\g_1=(\h_1(y,x))^2_x$, where $\h_1$ is the branch from $(y,x)$ in the subtree generated by the ball $B_R(x)$ in $T_{D+1}:=(V^{D+1},E^{D+1})$ with $(y,x)\in E^{D+1}$.\label{lar tree1}
\item If $L=2R+1$ for some positive integer $R$, then
\begin{align}\label{bigtree2}
\lambda_2(\g)\geq \dfrac{2}{2\sum\limits_{i=0}^{R-1}D^i+D^R}=\dfrac{2(D-1)}{D^{R+1}+D^R-2}.
\end{align}
Moreover, the equality holds if and only if $\g=\f:=(\tilde V,\tilde E)$, where $\f$ is exactly the subtree generated by $B_R(z)\bigcup B_R(w)$ in $T_{D+1}:=(V^{D+1},E^{D+1})$ with $(z,w)\in E^{D+1}$.\label{lar tree2}
\end{enumerate}

\tmd

\pf
For (\ref{lar tree1}),
one easily shows that such $\g$ can be viewed as a subgraph of the graph $\g_2$ generated by $B_R(x)$ in $T_{D+1}$. By Theorem \ref{sigmasteklov}, direct calculus shows that
\begin{align}\label{bigtree3}
\lambda_2(\g_1)=\lambda_2(\g_2)=\sigma(\h_1(y,x),x)=\dfrac{1}{\sum\limits_{i=0}^{R-1}D^i}=\dfrac{D-1}{D^R-1}.
\end{align}
Thus (\ref{bigtree1}) follows from Theorem \ref{main1} and (\ref{bigtree3}).

We proceed to discuss the remaining claim.

($\Longleftarrow$) Applying Theorem \ref{main1}, we have $\lambda_2(\g_2)\leq\lambda_2(\g)\leq\lambda_2(\g_1)$. Hence we finish the proof by (\ref{bigtree3}).

($\Longrightarrow$) Let $d(x,x_j)=1$ with $1\leq j\leq \deg(x)$ for $\g$. Note that $\g^j(x_j,x)$ can be viewed as a subgraph of $\h_1(y,x)$, where $\g^j$ is the branch from $(x_j,x)$ in $\g$. By Corollary \ref{smono}, we have $\sigma(\g^j(x_j,x),x)\geq\sigma(\h_1(y,x),x)$. Then by Theorem \ref{main1}, Theorem \ref{sigmasteklov}, we get

$$\lambda_2(\g)\geq\lambda_2((\g)^2_x)=\min\limits_{1\leq j\leq\deg(x)}\sigma(\g^j(x_j,x),x)\geq\sigma(\h_1(y,x),x).$$
On the other hand, it is obvious that $\lambda_2(\g)=\dfrac{D-1}{D^R-1}=\sigma(\h_1(y,x),x)$ by (\ref{bigtree3}). Using Theorem \ref{sigmasteklov}, there exist $s,t$ such that $\sigma(\g^s(x_s,x),x)=\sigma(\g^t(x_t,x),x)=\sigma(\h_1(y,x),x)$ with $1\leq s,t\leq\deg(x)$. Applying Corollary \ref{smono}, we conclude the proof.

For (\ref{lar tree2}), it is similar to the proof of (\ref{lar tree1}). Note that $\lambda_2(\f)=\dfrac{2(D-1)}{D^{R+1}+D^R-2}$. It is clear that $\g$ can be viewed as a subgraph of $\f$, hence we deduce (\ref{bigtree2}). So we only need to show the equality case.

($\Longleftarrow$) This direction is trivial.

($\Longrightarrow$)
Recall that $\g\subset\f.$ A direct computation shows that $\sigma(\f_1(z,w),w)=\sigma(\f_2(w,z),z)=\dfrac{D-1}{D^{R+1}-1}<\dfrac{2(D-1)}{D^{R+1}+D^R-2}=\lambda_2(\f)$, then the proof follows from Corollary \ref{strongrigidity}.

\pfd

%%%%%%%%%%%%%%%%%%%%%%%%%%%%%
%%%%%%%%%%%%%%%%%%%%%%%%%%%%%
%%%%%%%%%%%%%%%%%%%%%%%%%%%%%

\bibliography{py2}
\bibliographystyle{alpha}

\end{document}